\documentclass[11pt]{article}
\usepackage{amsmath, amssymb, amsthm, fullpage}
\date{}
\title{\vspace{-1cm} Large subgraphs without complete bipartite graphs
}
\author{David Conlon\thanks{
Mathematical Institute, Oxford OX2 6GG, UK. Email:
{\tt david.conlon@maths.ox.ac.uk}.} \and Jacob Fox
\thanks{
Department of Mathematics, MIT, Cambridge, MA 02139-4307. Email:
{\tt fox@math.mit.edu}.} 
\and
Benny Sudakov \thanks{Department of Mathematics, ETH, 8092 Zurich, Switzerland and Department of Mathematics, UCLA, Los Angeles, CA 90095.
Email: {\tt benjamin.sudakov@math.ethz.ch}. 
}
}

\theoremstyle{plain}
\newtheorem{THM}{Theorem}[section]
\newtheorem*{THM*}{Theorem}
\newtheorem{PROP}[THM]{Proposition}
\newtheorem{LEMMA}[THM]{Lemma}

\theoremstyle{definition}

\begin{document}
\maketitle

\begin{abstract}
In this note, we answer the following question of Foucaud, Krivelevich and Perarnau.
What is the size of the largest $K_{r,s}$-free subgraph one can guarantee in every graph $G$ with $m$ edges?
We also discuss the analogous problem for hypergraphs.
\end{abstract}

\section{Introduction}
Motivated by the classical Tur\'an problem, Foucaud, Krivelevich and Perarnau \cite{FKP} proposed to study the size of the largest $H$-free
subgraph one can always find in every graph $G$ with $m$ edges. Denote this function by $f(m,H)$. It is easy to determine $f(m,H)$ asymptotically if $H$ is not bipartite. In \cite{FKP}, the authors studied this problem when forbidding all even cycles in the subgraph up to length $2k$ and obtained estimates that are tight up to a logarithmic factor. They also asked to determine $f(m,H)$ when $H$ is a complete bipartite graph. The goal of this note is to resolve this question. 

\section{Complete bipartite graphs}
Let $K_{r,s}$ be the complete bipartite graph with parts of order $r$ and $s$, where $2 \leq r \leq s$. The following theorem gives a lower bound on $f(m,K_{r,s})$.
\begin{THM}
\label{th1}
Every graph $G$ with $m$ edges contains a $K_{r,r}$-free subgraph of size at least $\frac{1}{4}m^{\frac{r}{r+1}} $.
\end{THM}

To prove this theorem we need un upper bound on the maximum number of copies of $K_{r,r}$ which one can find in a graph with $m$ edges.
The problem of maximizing the number of copies of a fixed graph $H$ was solved by Alon \cite{Al} for all graphs and by
Friedgut and Kahn \cite{FK} for all hypergraphs. For our purposes the following easy estimate will suffice.

\begin{LEMMA}
\label{l1}
Every graph $G$ with $m$ edges contains at most $2m^r$ copies of $K_{r,r}$.
\end{LEMMA}

\noindent
{\bf Proof.}\, Note that every copy of $K_{r,r}$ in $G$ contains a matching of size $r$. Clearly the number of such matchings in $G$ is at most 
${m \choose r}$. Also note that every matching in $G$ of size $r$ can appear in at most $2^r$ copies of $K_{r,r}$. This implies that the total number of such copies is at most $2^r{m \choose r} \leq 2m^r$. \hfill $\Box$ \vspace{0.1cm}

\noindent
Using this lemma, together with a simple probabilistic argument, one can prove a lower bound on $f(m,K_{r,s})$.
\vspace{0.1cm}

\noindent{\bf Proof of Theorem \ref{th1}.}\, Let $G$ be a graph with $m$ edges. Consider a random subgraph $G'$ of $G$, obtained by choosing
every edge randomly and independently with probability $p=\frac{1}{2}m^{-1/(r+1)}$.
Then the expected number of edges in $G'$ is $mp$. Also, by
Lemma \ref{l1}, the expected number of copies of $K_{r,r}$ in $G'$ is at most $2p^{r^2}m^{r}$. Delete one edge from every copy of 
$K_{r,r}$ contained in $G'$. This gives a $K_{r,r}$-free subgraph of $G$, which by linearity of expectation, has at least
$$pm-2p^{r^2}m^{r} \geq \frac{1}{2} m^{\frac{r}{r+1}}-\frac{1}{8} m^{\frac{r}{r+1}} \geq \frac{1}{4} m^{\frac{r}{r+1}}$$
edges on average. Hence, there exists a choice of $G'$ which produces a $K_{r,r}$-free subgraph of $G$ of size at least
$\frac{1}{4}m^{\frac{r}{r+1}} $. \hfill $\Box$
\vspace{0.1cm}

Next we show that this gives an estimate on $f(m,K_{r,s})$ which is tight up to a constant factor depending on $s$ by taking $G$ to be an appropriately chosen complete bipartite graph with $m$ edges.
\begin{THM}
\label{th2}
Let $2 \leq r \leq s$ and let $G$ be a complete bipartite graph with parts $U$ and $W$, where $|U|=m^{1/(r+1)}$ and $|W|=m^{r/(r+1)}$. Then $G$ has $m$ edges and the largest $K_{r,s}$-free subgraph of $G$ has at most $s m^{r/(r+1)}$ edges.
\end{THM}

\noindent{\bf Proof.}\, The proof is a simple application of the counting argument of 
K\H{o}v\'ari-S\'os-Tur\'an \cite{KST}. Let $G'$ be a $K_{r,s}$-free subgraph of $G$ and let $d=e(G')/|W|$ be the average degree of vertices of $G'$ in $W$. If $d \geq s$, then, by convexity,
$$\sum_{w \in W} { d_{G'}(w) \choose r} \geq |W| {d \choose r} \geq {s \choose r}m^{r/(r+1)}\geq s m^{r/(r+1)}/r!\,.$$
On the other hand, since $G'$ is $K_{r,s}$-free we have that
$$\sum_{w \in W} {d_{G'}(w) \choose r} < s {|U| \choose r} \leq s |U|^r/r!=s m^{r/(r+1)}/r!\,.$$
This contradiction completes the proof of the theorem.  \hfill $\Box$
\vspace{0.1cm}

\noindent{\bf Remarks.}
\begin{itemize}
\item 
Since $K_{2,2}$ is also a $4$-cycle, our result improves by a logarithmic factor an estimate obtained by
Foucaud, Krivelevich and Perarnau \cite{FKP}.
\item
Since the Tur\'an number for $K_{r,s}$ is not known in general,
it is somewhat surprising that one can prove a tight bound on the size of the largest $K_{r,s}$-free subgraph in graphs with $m$ edges. \end{itemize}

\section{Hypergraphs}
The results presented in the previous section can be extended to $k$-uniform hypergraphs, which, for brevity, we call $k$-graphs.
Given a fixed $k$-graph $H$, let $f(m,H)$ denote the size of the largest $H$-free
subgraph one can always find in every $k$-graph $G$ with $m$ edges. Let $K^{(k)}_{r,\ldots,r}$ denote the complete $k$-partite
$k$-graph with parts of size $r$.

 \begin{THM}
 \label{th3}
 Every $k$-graph $G$ with $m$ edges contains a $K^{(k)}_{r,\ldots,r}$-free subgraph of size at least $\frac{1}{4}m^{\frac{q-1}{q}} $, where
 $q=\frac{r^k-1}{r-1}$.
 \end{THM}

\noindent{\bf Proof.}\, Let $G$ be a $k$-graph with $m$ edges. Every copy of $K^{(k)}_{r,\ldots,r}$ in $G$ contains a matching of size $r$ and the number of such matchings is at most ${m \choose r}$. On the other hand, every matching in $G$ of size $r$ can appear in at most $(k!)^r$ copies of $K_{r,r}$. This implies that the total number of such copies is at most $(k!)^r{m \choose r}$. 

Consider a random subgraph $G'$ of $G$, obtained by choosing
every edge randomly and independently with probability $p=\frac{1}{2}m^{-1/q}$.
Then the expected number of edges in $G'$ is $mp$ and the expected number of copies of $K^{(k)}_{r,\ldots,r}$ in $G'$ is at most 
$(k!)^rp^{r^k}{m \choose r}$. Delete one edge from every copy of 
$K^{(k)}_{r,\ldots,r}$ contained in $G'$. This gives a $K^{(k)}_{r,\ldots,r}$-free subgraph of $G$ with at least
$$pm- (k!)^rp^{r^k}{m \choose r} \geq \frac{1}{4}m^{\frac{q-1}{q}} $$
expected edges. Hence, there exists a choice of $G'$ which produces a $K^{(k)}_{r,\ldots,r}$-free subgraph of $G$ of this size.
\hfill $\Box$ \vspace{0.1cm}

We can again see that this estimate is tight up to a constant factor depending on $r$.

\begin{THM}
\label{th4}
Let $2 \leq r, k$,  $q=\frac{r^k-1}{r-1}$ and let $G$ be a complete $k$-partite $k$-graph with parts $U_i, 1\leq i \leq k$, such that $|U_i|=m^{r^{i-1}/q}$. Then $G$ has $m$ edges and the largest $K^{(k)}_{r,\ldots,r}$-free subgraph of $G$ has at most $rm^{(q-1)/q}$ edges.
\end{THM}

The proof of this theorem uses a similar counting argument to the graph case but is more involved. 
It follows from the following statement, which one can prove by induction. This technique has its origins in a paper of Erd\H{o}s \cite{E64}.

\begin{PROP}
Let $G$ be a $k$-partite $k$-graph with parts $U_i, 1\leq i \leq k$, such that $|U_i|=n^{r^i}$ and with $a \prod_{i \geq 2}|U_i|$ edges and $a \geq r$. Then $G$ contains at least ${a \choose r}\prod_{i \leq k-1} {|U_i| \choose r}$ copies of $K^{(k)}_{r,\ldots,r}$.
\end{PROP}

\noindent{\bf Proof.}\, We prove this by induction on $k$. The base case $k=1$ is trivial, by properly interpreting empty products as one. 

Now suppose we know the statement for $k-1$. For every vertex $x \in U_k$, denote by $G_x$ the $(k-1)$-partite $(k-1)$-graph which is the link of vertex $x$ (i.e., the collection of all subsets of size $k-1$ which together with $x$ form an edge of $G$). Let 
$a_x \prod_{i=2}^{k-1}|U_i|$ be the number of edges in $G_x$. By definition, $\sum_x a_x=a|U_k|=an^{r^k}$. By the induction hypothesis, each 
$G_x$ contains at least ${a_x \choose r}\prod_{i \leq k-2} {|U_i| \choose r}$ copies of 
$K^{(k-1)}_{r,\ldots,r}$. By convexity, the total number of such copies added over all $G_x$ is at least
$${a \choose r}n^{r^k}\prod_{i \leq k-2} {|U_i| \choose r} ={a \choose r}|U_{k-1}|^r\prod_{i \leq k-2} {|U_i| \choose r}
\geq r!{a \choose r} \prod_{i \leq k-1} {|U_i| \choose r} \geq
a\prod_{i \leq k-1} {|U_i| \choose r}.$$
For every subset $S$ which intersects each $U_i$ with $i \leq k-1$ in exactly $r$ vertices, denote by $d(S)$ the number of vertices $x \in U_k$ such that
$x$ forms an edge of $G$ together with every subset of $S$ of size $k-1$ which contain one vertex from every $U_i$. By the above discussion, we have that $\sum_Sd(S) \geq  a\prod_{i \leq k-1} {|U_i| \choose r}$, that is, at least the number of all copies of $K^{(k-1)}_{r,\ldots,r}$ in all $G_x$.
On the other hand, by the definition of $d(S)$,  the number of copies of $K^{(k)}_{r,\ldots,r}$ in $G$ equals
$\sum_S {d(S) \choose r}$. Since the total number of sets $S$ is $\prod_{i \leq k-1} {|U_i| \choose r}$, the result now follows by convexity.
\hfill $\Box$

\vspace{0.4cm}
\noindent
{\bf Acknowledgment.}\, We would like to thank M. Krivelevich for bringing this problem to our attention and for sharing with us preprint 
\cite{FKP}.


\begin{thebibliography}{99}

\bibitem{Al}
N. Alon, On the number of subgraphs of prescribed type of graphs with a given number of edges, {\em Israel J. Math.} {\bf 38} (1981), 116--130. 
\bibitem{E64}
P. Erd\H{o}s, On extremal problems of graphs and generalized graphs, {\em Israel J. Math.} {\bf 2} (1964), 183--190.
\bibitem{FKP}
F. Foucaud, M. Krivelevich and G. Perarnau, Large subgraphs without short cycles, preprint.
\bibitem{FK}
E. Freidgut and J. Kahn, On the number of copies of one hypergraph in another,
{\em Israel J. Math.} {\bf 105} (1998), 251--256. 
\bibitem{KST}
T. K\H{o}v\'ari, V.T. S\'os and P. Tur\'an, On a problem of K. Zarankiewicz,
{\em Colloq. Math.} {\bf 3} (1954), 50--57. 
\end{thebibliography}
\end{document}